\documentclass[11pt]{article}
\usepackage{amsfonts}
\usepackage{bbm}
\usepackage{mathrsfs}
\usepackage{amscd}
\usepackage{color}
\usepackage{amsmath,amsfonts,amssymb,amscd}
\usepackage{indentfirst,graphics,epsfig,psfrag}
\input{epsf}
\usepackage{ifpdf}
\usepackage{enumerate}
\usepackage{appendix}
\usepackage{enumerate}

\usepackage{lineno}

\makeatletter \@addtoreset{figure}{section} \makeatother
\makeatletter
\long\def\@makecaption#1#2{%
   \vskip 10\p@
   \setbox\@tempboxa\hbox{{#1}\ \ #2}%
   \ifdim \wd\@tempboxa >\hsize

       {#1}\ \ #2\par
   \else
       \hbox to\hsize{\hfil\box\@tempboxa\hfil}%
   \fi}
\makeatother

\begin{document}
\title{\textbf{Resistance distance and Kirchhoff index
in generalized $R$-vertex and $R$-edge corona for graphs}}
\author{
\small  Qun Liu~$^{a}$\\
\small  a. School of Mathematics and Statistics, \\
\small Hexi University, Gansu, Zhangye, 734000, P.R. China\\
\small E-mails: Liuqun@fudan.edu.cn}
\date{}
\maketitle
\begin{abstract}
For a graph $G$, the
graph $R(G)$ of a graph $G$ is the graph obtained by adding a new vertex
for each edge of $G$ and joining each new vertex to
both end vertices of the corresponding edge.
Let $I(G)$ be the set of newly added vertices, $i.e$ $I(G)=V(R(G))\setminus V(G)$.
The generalized $R$-vertex corona of $G$ and
$H_{i}$ for $i=1,2,...,n$, denoted by $R(G)\boxdot \wedge_{i=1}^{n}H_{i}$,
is the graph obtained from $R(G)$ and $H_{i}$ by joining
the ith vertex of $V(G)$ to every vertex in $H_{i}$.
The generalized $R$-edge corona of $G$ and
$H_{i}$ for $i=1,2,...,m$, denoted by $R(G)\ominus\wedge_{i=1}^{m}H_{i}$,
is the graph obtained from $R(G)$ and $H_{i}$ by joining
the ith vertex of $I(G)$ to every vertex in $H_{i}$.
In this paper, we derive closed-form formulas for resistance distance
and Kirchhoff index of $R(G)\boxdot \wedge_{i=1}^{n}H_{i}$
and $R(G)\ominus\wedge_{i=1}^{m}H_{i}$ whenever $G$ and
$H_{i}$ are arbitrary graph. These results
generalize the existing results in $\cite{HPZ}$.
\\[2mm]

\\[1mm]

\end{abstract}

\section{ Introduction }
All graphs considered in this paper are simple and undirected.
The resistance distance between vertices $u$ and $v$ of $G$ was defined
by Klein and Randi$\acute{c}$ $\cite{KR}$ to be the effective resistance between nodes $u$ and $v$ as computed with Ohm's law when all the edges of $G$
are considered to be unit resistors.
The Kirchhoff index $Kf(G)$ was defined in $\cite{KR}$
as $Kf(G)=\sum_{u<v}r_{uv}$, where $r_{uv}(G)$ denotes the resistance distance between $u$ and $v$ in $G$. These novel parameters are in fact intrinsic to
the graph theory and has some nice properties and applications in chemistry.
For the study of resistance distance and Kirchhoff index, one may be referred to
the recent works
$(\cite{BG}, \cite{BRB}, \cite{BYZhZH}), \cite{ChZh}-\cite{YK})$ and the references therein.

Let $G=(V(G),E(G))$ be a graph with vertex set $V(G)$ and edge set $E(G)$. Let $d_{i}$ be the degree of vertex $i$ in $G$ and $D_{G}=diag(d_{1}, d_{2},\cdots d_{|V(G)|})$ the diagonal matrix with all vertex degrees of $G$ as its diagonal entries. For a graph $G$, let $A_{G}$ and $B_{G}$ denote the adjacency matrix and vertex-edge incidence matrix of $G$, respectively. The matrix $L_{G}=D_{G}-A_{G}$ is called the Laplacian matrix of $G$, where $D_{G}$ is the diagonal matrix of vertex degrees of $G$. We use $\mu_{1}(G)\geq u_{2}(G)\geq\cdots\geq\mu_{n}(G)=0$ to denote the eigenvalues of $L_{G}$.
For other undefined notations and terminology from graph theory, the readers may refer to $\cite{RB}$ and the references therein.

In $\cite{PL}$, Lu et.al generalize the
corona operation and define the generalized $R$-vertex corona. For a graph $G$, the
graph $R(G)$ of a graph $G$ is the graph obtained by adding a new vertex
for each edge of $G$ and joining each new vertex to
both end vertices of the corresponding edge.
Let $I(G)$ be the set of newly added vertices, $i.e$ $I(G)=V(R(G))\setminus V(G)$.

{\bf Definition $1.1$}$(\cite{PL})$
The generalized $R$-vertex corona of $G$ and
$H_{i}$ for $i=1,2,...,n$, denoted by $R(G)\boxdot \wedge_{i=1}^{n}H_{i}$,
is the graph obtained from $R(G)$ and $H_{i}$ by joining
the ith vertex of $V(G)$ to every vertex in $H_{i}$.

{\bf Definition $1.2$}
The generalized $R$-edge corona of $G$ and
$H_{i}$ for $i=1,2,...,m$, denoted by $R(G)\ominus\wedge_{i=1}^{m}H_{i}$,
is the graph obtained from $R(G)$ and $H_{i}$ by joining
the ith vertex of $I(G)$ to every vertex in $H_{i}$.

Bu et al. investigated resistance distance in
subdivision-vertex join and subdivision-edge join
of graphs $\cite{BYZhZH}$. Liu et al. $\cite{LZhB}$ gave the resistance
distance and Kirchhoff index of $R$-vertex
join and $R$-edge join of two graphs. In $\cite{LP}$,
the resistance distance of subdivision-vertex and subdivision-edge
coronae are obtained. Motivated by the results, in this paper we
considered the generalization of the $R$-vertex corona and
the $R$-edge corona to the case of $n$$(m)$ different graphs
and we obtain the resistances distance and the Kirchhoff index in
terms of the corresponding parameters of the factors. These results
generalize the existing results
in $\cite{HPZ}$.

\section{Preliminaries}
The $\{1\}$-inverse of $M$ is a matrix $X$ such that $MXM=M$. If $M$ is singular, then it has infinite
$\{1\}$- inverse $\cite{A. Ben-Israel}$.
For a square matrix $M$, the group inverse of $M$, denoted by $M^{\#}$, is the unique matrix $X$ such that $MXM=M$, $XMX=X$ and $MX=XM$. It is known that $M^{\#}$ exists if and only if $rank(M)=rank(M^{2})$ $(\cite{A. Ben-Israel},\cite{BSZHW})$. If $M$ is real symmetric, then $M^{\#}$ exists and $M^{\#}$ is a symmetric $\{1\}$- inverse of $M$. Actually, $M^{\#}$ is equal to the Moore-Penrose inverse of $M$ since $M$ is symmetric $\cite{BSZHW}$.

It is known that resistance distances in a connected graph $G$ can be obtained from any $\{1\}$- inverse of $G$ $(\cite{BG})$. We use $M^{(1)}$ to denote any $\{{1}\}$- inverse of a matrix $M$, and let $(M)_{uv}$ denote the $(u,v)$- entry of $M$.

{\bf Lemma 2.1} $(\cite{BSZHW})$\ Let $G$ be a connected graph. Then
$$
r_{uv}(G)=(L^{(1)}_{G})_{uu}+(L^{(1)}_{G})_{vv}-(L^{(1)}_{G})_{uv}
-(L^{(1)}_{G})_{vu}=(L^{\#}_{G})_{uu}+(L^{\#}_{G})_{vv}-2(L^{\#}_{G})_{uv}.$$

Let $\bf{1_{n}}$ denotes the column vector of dimension $n$ with all the entries equal one. We will often use $\bf{1}$ to denote an all-ones column vector if the dimension can be read from the context.

{\bf Lemma 2.2} $(\cite{BYZhZH})$ \ For any graph $G$, we have
$L^{\#}_{G}\bf{1}$$=0.$

{\bf Lemma 2.3} $(\cite{FZZh})$ \ Let
 \[
\begin{array}{crl}
M=\Bigg(
  \begin{array}{cccccccccccccccc}
   A& B  \\
   C & D \\
 \end{array}
  \Bigg)
\end{array}
\]
be a nonsingular matrix. If $A$ and $D$ are nonsingular,
then
\[
\begin{array}{crl}
M^{-1}&=&\Bigg(
  \begin{array}{cccccccccccccccc}
   A^{-1}+A^{-1}BS^{-1}CA^{-1}&-A^{-1}BS^{-1} \\
   -S^{-1}CA^{-1} & S^{-1}\\
 \end{array}
  \Bigg)
\\&=&\Bigg(
  \begin{array}{cccccccccccccccc}
   (A-BD^{-1}C)^{-1} &-A^{-1}BS^{-1} \\
   -S^{-1}CA^{-1} & S^{-1}\\
 \end{array}
  \Bigg),
\end{array}
\]
where $S=D-CA^{-1}B.$

For a square matrix $M$, let $tr(M)$ denote the trace of $M$.

{\bf Lemma 2.4} $(\cite{SWZhB})$ \ Let $G$ be a connected graph on $n$ vertices. Then
$$Kf(G)=ntr(L^{(1)}_{G})-1^{T}L^{(1)}_{G}1=ntr(L^{\#}_{G}).$$

{\bf Lemma 2.5}$(\cite{DJ})$ \ Let $G$ be a connected graph of order $n$ with edge set $E$. Then
$$\sum_{u<v, uv\in E}r_{uv}(G)=n-1.$$
For a vertex $i$ of a graph $G$, let $T(i)$ denote the set of all neighbors of $i$ in $G$.

{\bf Lemma 2.6}$(\cite{BYZhZH})$ \ Let $G$ be a connected graph. For any $i,j\in V(G)$,
$$r_{ij}(G)=d^{-1}_{i}(1+\sum_{k\in T(i)}r_{kj}(G)-d^{-1}_{i}\sum_{k,l\in T(i)}r_{kl}(G)).$$

{\bf Lemma 2.7} $(\cite{LZhB})$ \ Let $G$ be a graph of order
$n$. For any $a,b>0$ satisfying $b\neq a$, we have
$$(L_{G}+aI_{n}-\frac{a}{b}j_{n\times n})^{-1}
=(L_{G}+aI_{n})^{-1}+\frac{1}{a(b-n)}j_{n\times n},$$
where $j_{n\times n}$ denotes the $n\times n$ matrix
with all entries equal to one.

{\bf Lemma 2.8} $(\cite{QL})$  \ Let
 \[
\begin{array}{crl}
L=\left(
  \begin{array}{cccccccccccccccc}
   A& B \\
   B^{T} & D \\
 \end{array}
  \right)
\end{array}
\]
be a symmetric block matrix.
If $D$ is nonsingular, then
 \[
\begin{array}{crl}
X=\left(
  \begin{array}{cccccccccccccccc}
   H^{\#}& -H^{\#}BD^{-1} \\
   -D^{-1}B^{T}H^{\#} & D^{-1}+D^{-1}B^{T}H^{\#}BD^{-1} \\
 \end{array}
  \right)
\end{array}
\]
is a symmetric $\{1\}$-inverse of $L$,
where $H=A-BD^{-1}B^{T}$.

{\bf Lemma 2.9} $(\cite{KR})$ \ Let $k$ be a cut-vertex
of a graph, and let $i$ and $j$ be vertices occurring in
different components which arise upon deletion of $k$. Then
$$r_{ij}=r_{ik}+r_{kj}.$$

\section{The resistance distance and Kirchhoff index of
$R(G)\boxdot \wedge_{i=1}^{n}H_{i}$}
In this section, we focus on determing the resistance distance and Kirchhoff index of generalized $R$-vertex corona
$R(G)\boxdot \wedge_{i=1}^{n}H_{i}$ whenever $G$ and $H_{i}$$(i=1,2,...,n)$ be an arbitrary graph.

{\bf Theorem 3.1} \ Let $G$ be a connected graph with $n$ vertices and $m$
edges, Let $H_{i}$ be a graph with $t_{i}$ vertices for
$i=1,2,...,n$.
Then $R(G)\boxdot \wedge_{i=1}^{n}H_{i}$ have the resistance distance
and Kirchhoff index as follows:

(i) For any $i,j\in V(G)$, we have
\begin{eqnarray*}
r_{ij}(L_{R(G)\boxdot \wedge_{i=1}^{n}H_{i}})&=&\frac{2}{3}(L^{\#}_{G})_{ii}+\frac{2}{3}(L^{\#}_{G})_{jj}
-\frac{4}{3}(L^{\#}_{G})_{ij}=\frac{2}{3}r_{ij}(G),
\end{eqnarray*}

(ii) For any $i,j\in V(H_{k})$$(k=1,2,...,n)$, we have
\begin{eqnarray*}
r_{ij}(L_{R(G)\boxdot \wedge_{i=1}^{n}H_{i}})&=&((L_{H_{k}}+I_{t_{k}})^{-1})_{ii}+
((L_{H_{k}}+I_{t_{k}})^{-1})_{jj}
-2((L_{H_{k}}+I_{t_{k}})^{-1})_{ij}.
\end{eqnarray*}

(iii) For any $i,j\in R(G)$,
we have
\begin{eqnarray*}
r_{ij}(R(G)\boxdot \wedge_{i=1}^{n}H_{i})
=\frac{2}{3}r_{ij}(G).
\end{eqnarray*}

(iv) For any $i\in V(G)$, $j\in V(H_{k})$$(k=1,2,...,n)$, we have
\begin{eqnarray*}
r_{ij}(R(G)\boxdot \wedge_{i=1}^{n}H_{i})
&=&r_{ik}(R(G))+r_{kj}(F_{k}),
\end{eqnarray*}
where $F_{k}=H_{k}\vee \{v\}$.

(v) For any $i\in V(H_{k})$, $j\in V(H_{l})$, we have
\begin{eqnarray*}
r_{ij}(R(G)\boxdot \wedge_{i=1}^{n}H_{i})
&=&r_{kl}(R(G))+r_{ik}(F_{k})+r_{jl}(F_{l}),
\end{eqnarray*}
where $F_{k}=H_{k}\vee \{v\}$.

(vi)
$Kf(R(G)\boxdot \wedge_{i=1}^{n}H_{i})$
\begin{eqnarray*}
&=&(n+2m+\sum_{i=1}^{n}t_{i})\left(\frac{2}{3n}Kf(G)
+\frac{m}{2}+\frac{1}{2}tr(D_{G}L^{\#}_{G})-\frac{n-1}{4}+
\sum_{i=1}^{n}\sum_{j=1}^{t_{i}}\frac{1}{\mu_{i}(H_{j})+1}\right.\\
&&\left.+2tr(Q^{T}L^{\#}_{G}Q)\right)
-\left(\frac{m}{2}+\frac{1}{4}\pi^{T}L^{\#}_{G}\pi+\pi^{T}L^{\#}_{G}\delta
+\sum_{i=1}^{n}t_{i}+\delta^{T}L^{\#}_{G}\delta\right),
\end{eqnarray*}
where $Q$ equals (3.1), $\pi^{T}=(d_{1},d_{2},...,d_{n})$, $\delta^{T}=(t_{1},t_{2},...,t_{n})$.

{\bf Proof} \ Let $R(G)$ and $D_{G}$ be the incidence matrix
and degree matrix of $G$. With a suitable labeling for vertices of $R(G)\boxdot \wedge_{i=1}^{n}H_{i}$,
the Laplacian matrix of $R(G)\boxdot \wedge_{i=1}^{n}H_{i}$ can be written as follows:
\[
\begin{array}{crl}
 L_{R(G)\boxdot \wedge_{i=1}^{n}H_{i}}=\left(
  \begin{array}{cccccccccccccccc}
   P+L_{G}& -R(G)& -Q\\
  -R^{T}(G) & 2I_{m}& 0 \\
  -Q^{T}& 0 & T\\
 \end{array}
  \right),
\end{array}
\]
where
\begin{eqnarray}
 P=\left(
  \begin{array}{cccccccccccccccc}
   d_{1}+t_{1}& 0 & 0&...&0\\
  0&d_{2}+t_{2}&0&...&0\\
  0&0& ...&...&0\\
   0&0& 0&...&d_{n}+t_{n}\\
 \end{array}
  \right),~~~~~~Q=
  \left(
  \begin{array}{cccccccccccccccc}
   1^{T}_{t_{1}}& 0 & 0&...&0\\
  0&1^{T}_{t_{2}}&0&...&0\\
  0&0& ...&...&0\\
   0&0& 0&...&1^{T}_{t_{n}}\\
 \end{array}
  \right),
\end{eqnarray}

\[
\begin{array}{crl}
 T=\left(
  \begin{array}{cccccccccccccccc}
   L_{H_{1}}+I_{t_{1}}& 0 & 0&...&0\\
  0&L_{H_{2}}+I_{t_{2}}&0&...&0\\
  0&0& ...&...&0\\
   0&0& 0&...&L_{H_{n}}+I_{t_{n}}\\
 \end{array}
  \right).
\end{array}
\]
First we begin with the computation of $\{1\}$-inverse of $R(G)\boxdot \wedge_{i=1}^{n}H_{i}$.

By Lemma 2.8, we have
\[
\begin{array}{crl}
 H&=&L_{G}+P-\left(
  \begin{array}{cccccccccccccccc}
  -R(G)& -Q \\
 \end{array}
  \right)
 \left(
  \begin{array}{cccccccccccccccc}
   \frac{1}{2}I_{m}&0\\
   0&T^{-1}\\
 \end{array}
  \right)
  \left(
  \begin{array}{cccccccccccccccc}
    -R^{T}(G) \\
  -Q^{T}\\
 \end{array}
  \right)\\
  &=&L_{G}+P-\left(
  \begin{array}{cccccccccccccccc}
    -\frac{1}{2}R(G)&-QT^{-1}\\
 \end{array}
  \right)
  \left(
  \begin{array}{cccccccccccccccc}
     -R^{T}(G) \\
  -Q^{T}\\
 \end{array}
  \right)\\
  &=&L_{G}+D_{G}+\left(
  \begin{array}{cccccccccccccccc}
   t_{1}& 0 & 0&...&0\\
  0&t_{2}&0&...&0\\
  0&0& ...&...&0\\
   0&0& 0&...&t_{n}\\
 \end{array}
  \right)-\frac{1}{2}(D_{G}+A_{G})-\left(
  \begin{array}{cccccccccccccccc}
   t_{1}& 0 & 0&...&0\\
  0&t_{2}&0&...&0\\
  0&0& ...&...&0\\
   0&0& 0&...&t_{n}\\
 \end{array}
  \right)\\
 &=&\frac{3}{2}L_{G},
\end{array}
\]
so $H^{\#}=\frac{2}{3}L^{\#}_{G}$.

According to Lemma 2.8, we calculate $-H^{\#}BD^{-1}$ and $-D^{-1}B^{T}H^{\#}$.
\[
\begin{array}{crl}
-H^{\#}BD^{-1}&=&-\frac{2}{3}L^{\#}_{G}\left(
  \begin{array}{cccccccccccccccc}
 -R(G)& -Q \\
 \end{array}
  \right)
  \left(
  \begin{array}{cccccccccccccccc}
   \frac{1}{2}I_{m}&0\\
   0&T^{-1}\\
 \end{array}
  \right)\\
  &=&-\frac{2}{3}L^{\#}_{G}\left(
  \begin{array}{cccccccccccccccc}
 -\frac{1}{2}R(G)&-QT^{-1}\\
 \end{array}
  \right)
   =\left(
  \begin{array}{cccccccccccccccc}
 \frac{1}{3}L^{\#}_{G}R(G)&\frac{2}{3}L^{\#}_{G}Q\\
 \end{array}
  \right)\\
\end{array}
\]
and
\[
\begin{array}{crl}
-D^{-1}B^{T}H^{\#}=-(H^{\#}BD^{-1})^{T}
  =\left(
  \begin{array}{cccccccccccccccc}
   \frac{1}{3}R^{T}(G)L^{\#}_{G}\\
   \frac{2}{3}Q^{T}L^{\#}_{G} \\
 \end{array}
  \right).\\
\end{array}
\]
We are ready to compute the $D^{-1}B^{T}H^{\#}BD^{-1}$.
\[
\begin{array}{crl}
D^{-1}B^{T}H^{\#}BD^{-1}
&=&\begin{array}{crl}
 \left(
  \begin{array}{cccccccccccccccc}
     \frac{1}{2}I_{m}&0\\
   0&T^{-1}\\
 \end{array}
 \right) \left(
  \begin{array}{cccccccccccccccc}
 -R^{T}(G)\\
  -Q^{T}\\
 \end{array}
 \right)L^{\#}_{G}
 \left(
  \begin{array}{cccccccccccccccc}
   -R(G)& -Q \\
 \end{array}
 \right)
 \left(
  \begin{array}{cccccccccccccccc}
     \frac{1}{2}I_{m}&0\\
   0&T^{-1}\\
 \end{array}
 \right)
\end{array}\\
&=&\left(
\begin{array}{cccccccccccccccc}
     \frac{1}{4}R^{T}(G)L^{\#}_{G}R(G)&\frac{1}{2}R^{T}(G)L^{\#}_{G}Q\\
   \frac{1}{2}Q^{T}L^{\#}_{G}R(G)&Q^{T}L^{\#}_{G}Q\\
 \end{array}
 \right).
\end{array}\\
\]
Based on Lemma 2.8, the following matrix
\begin{eqnarray}
N=\left(
  \begin{array}{cccccccccccccccc}
  \frac{2}{3}L^{\#}_{G}&\frac{1}{3}L^{\#}_{G}R(G)&\frac{2}{3}L^{\#}_{G}Q\\
   \frac{1}{3}R^{T}(G)L^{\#}_{G}&
   \frac{1}{2}I_{m}+\frac{1}{4}R^{T}(G)L^{\#}_{G}R(G)&
   \frac{1}{2}R^{T}(G)L^{\#}_{G}Q\\
   \frac{2}{3}Q^{T}L^{\#}_{G}&\frac{1}{2}Q^{T}L^{\#}_{G}R(G)&T^{-1}+Q^{T}L^{\#}_{G}Q\\
 \end{array}
  \right)
\end{eqnarray}
is a symmetric $\{1\}$- inverse of $L_{R(G)\boxdot \wedge_{i=1}^{n}H_{i}}$.

For any $i,j\in V(G)$, by Lemma 2.1 and the Equation $(3.2)$, we have
\begin{eqnarray*}
r_{ij}(L_{R(G)\boxdot \wedge_{i=1}^{n}H_{i}})&=&\frac{2}{3}(L^{\#}_{G})_{ii}+\frac{2}{3}(L^{\#}_{G})_{jj}
-\frac{4}{3}(L^{\#}_{G})_{ij}=\frac{2}{3}r_{ij}(G)
\end{eqnarray*}
as stated in $(i)$.

For any $i,j\in V(H_{k})$$(k=1,2,...,n)$, by Lemma 2.1 and the Equation $(3.2)$, we have
\begin{eqnarray*}
r_{ij}(L_{R(G)\boxdot \wedge_{i=1}^{n}H_{i}})&=&((L_{H_{k}}+I_{t_{k}})^{-1})_{ii}+
((L_{H_{k}}+I_{t_{k}})^{-1})_{jj}
-2((L_{H_{k}}+I_{t_{k}})^{-1})_{ij}
\end{eqnarray*}
as stated in (ii).

From the left side of above equation, we can obviously have
\begin{eqnarray*}
r_{ij}(F_{k})&=&((L_{H_{k}}+I_{t_{l}})^{-1})_{ii}+
((L_{H_{k}}+I_{t_{l}})^{-1})_{jj}
-2((L_{H_{k}}+I_{t_{l}})^{-1})_{ij},
\end{eqnarray*}
where $F_{k}=H_{k}\vee \{v\}$, i.e, $F_{k}$
is the graph obtained by
adding new edges from an isolated vetrtex
$v$ to every vertex of $H_{k}$.

For any $i,j\in R(G)$, by Lemma 2.1 and the Equation $(3.2)$,
we have
\begin{eqnarray*}
r_{ij}(R(G)\boxdot \wedge_{i=1}^{n}H_{i})
&=&r_{ij}(R(G)).
\end{eqnarray*}
By Lemma 3.1 in \cite{DC}, $r_{ij}(R(G))
=\frac{2}{3}r_{ij}(G)$, so $r_{ij}(R(G)\boxdot \wedge_{i=1}^{n}H_{i})
=\frac{2}{3}r_{ij}(G)$.

For any $i\in V(G)$, $j\in V(H_{k})$$(k=1,2,...,n)$, since
$i$ and $j$ belong to different components, then by Lemma 2.9,
we have
\begin{eqnarray*}
r_{ij}(R(G)\boxdot \wedge_{i=1}^{n}H_{i})
&=&r_{ik}(R(G))+r_{kj}(F_{k}).
\end{eqnarray*}

For any $i\in V(H_{k})$, $j\in V(H_{l})$,  by Lemma 2.9,
we have
\begin{eqnarray*}
r_{ij}(R(G)\boxdot \wedge_{i=1}^{n}H_{i})
&=&r_{kl}(R(G))+r_{ik}(F_{k})+r_{jl}(F_{l}).
\end{eqnarray*}

By Lemma 2.4, we have
\begin{eqnarray*}
Kf(L_{R(G)\boxdot \wedge_{i=1}^{nH_{i}}})
&=&(n+m+\sum_{i=1}^{n}t_{i})tr( N)-{\bf{1}^{T}}N{\bf{1}^{T}}\\
&=&(n+m+\sum_{i=1}^{n}t_{i})\left(\frac{2}{3}tr(L^{\#}_{G})
+tr\left(\frac{1}{2}I_{m}+\frac{1}{4}R^{T}(G)L^{\#}_{G}R(G)\right)+\right.\\
&&\left.+
tr(T^{-1}+Q^{T}L^{\#}_{G}Q)\right)-
{\bf{1}^{T}}N{\bf{1}^{T}}\\
&=&(n+m+\sum_{i=1}^{n}t_{i})\left(\frac{2}{3n}Kf(G)
+\frac{m}{2}+\frac{1}{4}\sum_{i<j,i,j\in E(G)}[(L^{\#}_{G})_{ii}+(L^{\#}_{G})_{jj}\right.\\
&&\left.+2(L^{\#}_{G})_{ij}]+tr\left(T^{-1}+Q^{T}L^{\#}_{G}Q\right)\right)
-{\bf{1}^{T}}N{\bf{1}^{T}}
\end{eqnarray*}
By Lemma 2.4, we get
\begin{eqnarray*}
Kf(L_{R(G)\boxdot \wedge_{i=1}^{n}H_{i}})
&=&(n+m+\sum_{i=1}^{n}t_{i})\left(\frac{2}{3n}Kf(G)
+\frac{m}{2}+\frac{1}{4}\sum_{i<j,i,j\in E(G)}[2(L^{\#}_{G})_{ii}+2(L^{\#}_{G})_{jj}\right.\\
&&\left.-r_{ij}(G)]+tr\left(T^{-1}+Q^{T}L^{\#}_{G}Q\right)\right)
-{\bf{1}^{T}}N{\bf{1}^{T}}\\
&=&(n+m+\sum_{i=1}^{n}t_{i})\left(\frac{2}{3n}Kf(G)
+\frac{m}{2}+\frac{1}{2}tr(D_{G}L^{\#}_{G})-\frac{n-1}{4}\right.\\
&&\left.+tr\left(T^{-1}+Q^{T}L^{\#}_{G}Q\right)\right)
-{\bf{1}^{T}}N{\bf{1}^{T}}
\end{eqnarray*}
Note that the eigenvalues of $(L(H_{i})+I_{t_{i}})$
$(i=1,2,...,n)$ are $\mu_{1}(H_{i})+1, \mu_{2}(H_{i})+1,...,\mu_{t_{i}}(H_{i})+1$. Then
\begin{eqnarray}
tr(T^{-1})=
\sum_{i=1}^{n}\sum_{j=1}^{t_{i}}\frac{1}{\mu_{i}(H_{j})+1}.
\end{eqnarray}
By Lemma 2.2, $L_{G}^{\#}\bf{1}$$=0$ and $({\bf{1}^{T}}\left(R^{T}(G)L^{\#}_{G}Q\right){\bf{1}})^{T}
={\bf{1}^{T}}\left(Q^{T}L^{\#}_{G}R(G)\right){\bf{1}}$, then
\begin{eqnarray*}
{\bf{1}^{T}}N{\bf{1}}&=&
\frac{m}{2}+
\frac{1}{4}{\bf{1}^{T}}\left(R^{T}(G)L^{\#}_{G}R(G)\right){\bf{1}}+
{\bf{1}^{T}}\left(R^{T}(G)L^{\#}_{G}Q\right){\bf{1}}
\\
&&+{\bf{1}^{T}}T^{-1}{\bf{1}}+{\bf{1}^{T}}\left(Q^{T}L^{\#}_{G}Q\right){\bf{1}}.
\end{eqnarray*}
Note that $R(G)\bf{1}=\pi$, where $\pi^{T}=(d_{1},d_{2},...,d_{n})$,
then ${\bf{1}^{T}}\left(R^{T}(G)L^{\#}_{G}R(G)\right){\bf{1}}
=\pi^{T}L^{\#}_{G}\pi$, so
\begin{eqnarray}
{\bf{1}^{T}}N{\bf{1}}&=&\frac{m}{2}+\frac{1}{4}\pi^{T}L^{\#}_{G}\pi+
\pi^{T}L^{\#}_{G}Q{\bf{1}}+{\bf{1}^{T}}T^{-1}{\bf{1}}+{\bf{1}^{T}}
\left(Q^{T}L^{\#}_{G}Q\right){\bf{1}}.
\end{eqnarray}
Let $R_{i}=L(H_{i})+I_{t_{i}}$$(i=1,2,...,n)$, then
\begin{eqnarray*}
{\bf{1}^{T}}T^{-1}{\bf{1}^{T}}&=&\left(
  \begin{array}{cccccccccccccccc}
 {\bf{1}^{T}_{t_{1}}}&{\bf{1}^{T}_{t_{2}}}&\cdots &{\bf{1}^{T}_{t_{n}}}\\
 \end{array}
  \right)
  \left(
  \begin{array}{cccccccccccccccc}
R_{1}^{-1}& 0 & 0&...&0\\
  0&R_{2}^{-1}&0&...&0\\
  0&0& ...&...&0\\
   0&0& 0&...&R_{n}^{-1}\\
 \end{array}
  \right)\left(
  \begin{array}{cccccccccccccccc}
   {\bf{1}_{t_{1}}}\\
  {\bf{1}_{t_{2}}}\\
  \cdots\\
  {\bf{1}_{t_{n}}}\\
 \end{array}
  \right)\\
  \end{eqnarray*}
  \begin{eqnarray}
 &=&\sum_{i=1}^{n}{\bf{1}^{T}_{t_{i}}}(L(H_{i})+I_{t_{i}})^{-1}{\bf{1}_{t_{i}}}
 =\sum_{i=1}^{n}t_{i},
\end{eqnarray}
and
\begin{eqnarray*}
{\bf{1}^{T}}Q^{T}&=&\left(
  \begin{array}{cccccccccccccccc}
 {\bf{1}^{T}_{t_{1}}}&{\bf{1}^{T}_{t_{2}}}&\cdots &{\bf{1}^{T}_{t_{n}}}\\
 \end{array}
  \right)
  \left(
  \begin{array}{cccccccccccccccc}
 {\bf{1}_{t_{1}}}& 0 & 0&...&0\\
  0& {\bf{1}_{t_{2}}}&0&...&0\\
  0&0& ...&...&0\\
   0&0& 0&...& {\bf{1}_{t_{n}}}\\
 \end{array}
  \right)\\
    \end{eqnarray*}
\begin{eqnarray}
 &=&(t_{1},t_{2},...,t_{n})=\delta^{T}.
\end{eqnarray}

Plugging $(3.3),(3.4),(3.5)$ and $(3.6)$ into $Kf(L_{R(G)\boxdot \wedge_{i=1}^{n}H_{i}})$,
we obtain the required result in $vi)$.

\section{The resistance distance and Kirchhoff index of
$R(G)\ominus\wedge_{i=1}^{m}H_{i}$}
In this section, we focus on determing the resistance distance and Kirchhoff index of generalized $R$-edge corona
$R(G)\ominus\wedge_{i=1}^{m}H_{i}$ whenever $G$ and $H_{i}$$(i=1,2,...,n)$ be an arbitrary graph.

{\bf Theorem 4.1} \ Let $G$ be a connected graph with $n$ vertices and $m$
edges, Let $H_{i}$ be a graph with $t_{i}$ vertices for
$i=1,2,...,m$.
Then $R(G)\ominus\wedge_{i=1}^{m}H_{i}$ have the resistance distance
and Kirchhoff index as follows:

(i) For any $i,j\in V(G)$, we have
\begin{eqnarray*}
r_{ij}(R(G)\ominus\wedge_{i=1}^{m}H_{i})
&=&\frac{2}{3}(L^{\#}_{G})_{ii}+\frac{2}{3}(L^{\#}_{G})_{jj}
-\frac{4}{3}(L^{\#}_{G})_{ij}=\frac{2}{3}r_{ij}(G),
\end{eqnarray*}

(ii) For any $i,j\in V(H_{k})$$(k=1,2,...,m)$, we have
\begin{eqnarray*}
r_{ij}(L_{R(G)\ominus\wedge_{i=1}^{m}H_{i}})&=&(L_{H_{k}}+I_{t_{k}}
-\frac{1}{2+t_{k}}j_{t_{k}})^{-1}_{ii}+
(L_{H_{k}}+I_{t_{k}}-\frac{1}{2+t_{k}}j_{t_{k}})^{-1}_{jj}\\&&
-2(L_{H_{1}}+I_{t_{k}}-\frac{1}{2+t_{k}}j_{t_{k}})^{-1}_{ij}.
\end{eqnarray*}

(iii) For any $i,j\in R(G)$,
we have
\begin{eqnarray*}
r_{ij}(R(G)\ominus\wedge_{i=1}^{m}H_{i})
=\frac{2}{3}r_{ij}(G).
\end{eqnarray*}

(iv) For any $i\in V(G)$, $j\in V(H_{k})$$(k=1,2,...,n)$, we have
\begin{eqnarray*}
r_{ij}(R(G)\ominus\wedge_{i=1}^{m}H_{i})
&=&r_{ik}(R(G))+r_{kj}(F_{k}),
\end{eqnarray*}
where $F_{k}=H_{k}\vee \{v\}$.

(v) For any $i\in V(H_{k})$, $j\in V(H_{l})$, we have
\begin{eqnarray*}
r_{ij}(R(G)\ominus\wedge_{i=1}^{m}H_{i})
&=&r_{kl}(R(G))+r_{ik}(F_{k})+r_{jl}(F_{l}),
\end{eqnarray*}
where $F_{k}=H_{k}\vee \{v\}$.

(vi)
$Kf(R(G)\ominus\wedge_{i=1}^{m}H_{i})$
\begin{eqnarray*}
&=&(n+2m+\sum_{i=1}^{n}t_{i})\left(\frac{2}{3n}Kf(G)
+\frac{m}{2}+\frac{1}{3}tr(D_{G}L^{\#}_{G})-\frac{n-1}{2}+\sum_{i=1}^{n}
\sum_{j=1}^{t_{i}}\frac{1}{\mu_{i}(H_{j})+1}\right.\\
&&\left.+\frac{2}{3}tr(F^{T}R^{T}(G)L^{\#}_{G}R(G)F)\right)
-\left(\frac{m}{2}+\frac{1}{6}\pi^{T}L^{\#}_{G}\pi+\sum_{i=1}^{m}t_{i}
+\frac{2}{3}\pi^{T}(G)L^{\#}_{G}R(G)\delta
+\pi^{T}L^{\#}_{G}\delta\right.\\&&\left.
+\frac{1}{2}\sum_{i=1}^{m}t_{i}(2+t_{i})+\frac{2}{3}
\delta^{T}R^{T}(G)L^{\#}_{G}R(G)\delta\right),
\end{eqnarray*}
where $F$ equals (4.7), $\pi^{T}=(d_{1},d_{2},...,d_{n})$, $\delta^{T}=(t_{1},0,...,0,t_{2},0,...,0...,t_{m})$.

{\bf Proof} \ Let $R(G)$ and $D_{G}$ be the incidence matrix
and degree matrix of $G$. With a suitable labeling for vertices of $R(G)\ominus\wedge_{i=1}^{m}H_{i}$,
the Laplacian matrix of $R(G)\ominus\wedge_{i=1}^{m}H_{i}$ can be written as follows:
\[
\begin{array}{crl}
 L_{R(G)\ominus\wedge_{i=1}^{m}H_{i}}=\left(
  \begin{array}{cccccccccccccccc}
   L_{G}+D_{G}& -R(G)& 0\\
  -R^{T}(G) & P& -M \\
  0& -M^{T} & Q\\
 \end{array}
  \right),
\end{array}
\]
where
\begin{eqnarray*}
 P=\left(
  \begin{array}{cccccccccccccccc}
   2+t_{1}& 0 & 0&...&0\\
  0&2+t_{2}&0&...&0\\
  0&0& ...&...&0\\
   0&0& 0&...&2+t_{m}\\
 \end{array}
  \right)_{m\times m},~~~~~~M=
  \left(
  \begin{array}{cccccccccccccccc}
   1^{T}_{t_{1}}& 0 & 0&...&0\\
  0&1^{T}_{t_{2}}&0&...&0\\
  0&0& ...&...&0\\
   0&0& 0&...&1^{T}_{t_{m}}\\
 \end{array}
  \right)_{m\times (t_{1}+t_{2}+\cdots +t_{n})},
\end{eqnarray*}

\begin{eqnarray*}
 Q=\left(
  \begin{array}{cccccccccccccccc}
   L_{H_{1}}+I_{t_{1}}& 0 & 0&...&0\\
  0&L_{H_{2}}+I_{t_{2}}&0&...&0\\
  0&0& ...&...&0\\
   0&0& 0&...&L_{H_{n}}+I_{t_{n}}\\
 \end{array}
  \right).
\end{eqnarray*}

Let
$A=L_{G}+D_{G}$,
 $B=\left(
  \begin{array}{cc}
-R(G)& 0\\
  \end{array}
\right)$,
 $B^{T}=\left(
  \begin{array}{cc}
   -R^{T}(G) \\
  0\\
  \end{array}
\right)$ and
$D=\left(
  \begin{array}{cccccccccccccccc}
  P&-M\\
  -M^{T}&Q\\
 \end{array}
  \right)$

First, we will compute $D^{-1}$.
By Lemma 2.3, we have
\[
\begin{array}{crl}
 S=\left(
  \begin{array}{cccccccccccccccc}
   L_{H_{1}}+I_{t_{1}}& 0 & 0&...&0\\
  0&L_{H_{2}}+I_{t_{2}}&0&...&0\\
  0&0& ...&...&0\\
   0&0& 0&...&L_{H_{m}}+I_{t_{m}}\\
 \end{array}
  \right)-\left(
  \begin{array}{cccccccccccccccc}
   1_{t_{1}}& 0 & 0&...&0\\
  0&1_{t_{2}}&0&...&0\\
  0&0& ...&...&0\\
   0&0& 0&...&1_{t_{m}}\\
 \end{array}
  \right)\\
  \left(
  \begin{array}{cccccccccccccccc}
   2+t_{1}& 0 & 0&...&0\\
  0&2+t_{2}&0&...&0\\
  0&0& ...&...&0\\
   0&0& 0&...&2+t_{n}\\
 \end{array}
  \right)^{-1}\left(
  \begin{array}{cccccccccccccccc}
   1^{T}_{t_{1}}& 0 & 0&...&0\\
  0&1^{T}_{t_{2}}&0&...&0\\
  0&0& ...&...&0\\
   0&0& 0&...&1^{T}_{t_{m}}\\
 \end{array}
  \right)
\end{array}
\]
\[
\begin{array}{crl}
S^{-1}=\left(
  \begin{array}{cccccccccccccccc}
   (L_{H_{1}}+I_{t_{1}}-\frac{1}{2+t_{1}}j_{t_{1}})^{-1}& 0 & 0&...&0\\
  0& (L_{H_{2}}+I_{t_{2}}-\frac{1}{2+t_{2}}j_{t_{2}})^{-1}&0&...&0\\
  0&0& ...&...&0\\
   0&0& 0&...& (L_{H_{m}}+I_{t_{m}}-\frac{1}{2+t_{m}}j_{t_{m}})^{-1}\
 \end{array}
  \right).
\end{array}
\]

According to Lemma 2.3, we have

\[
\begin{array}{crl}
 P-MQ^{-1}M^{T}=\left(
  \begin{array}{cccccccccccccccc}
   2+t_{1}& 0 & 0&...&0\\
  0&2+t_{2}&0&...&0\\
  0&0& ...&...&0\\
   0&0& 0&...&2+t_{m}\\
 \end{array}
  \right)-\left(
  \begin{array}{cccccccccccccccc}
   1^{T}_{t_{1}}& 0 & 0&...&0\\
  0&1^{T}_{t_{2}}&0&...&0\\
  0&0& ...&...&0\\
   0&0& 0&...&1^{T}_{t_{m}}\\
 \end{array}
  \right)\\
  \left(
  \begin{array}{cccccccccccccccc}
  (L_{H_{1}}+I_{t_{1}})^{-1}& 0 & 0&...&0\\
  0&(L_{H_{2}}+I_{t_{2}})^{-1}&0&...&0\\
  0&0& ...&...&0\\
   0&0& 0&...&(L_{H_{m}}+I_{t_{m}})^{-1}\\
 \end{array}
  \right)\left(
  \begin{array}{cccccccccccccccc}
   1_{t_{1}}& 0 & 0&...&0\\
  0&1_{t_{2}}&0&...&0\\
  0&0& ...&...&0\\
   0&0& 0&...&1_{t_{m}}\\
 \end{array}
  \right)\\
  =2I_{m},
\end{array}
\]
so $( P-MQ^{-1}M^{T})^{-1}=\frac{1}{2}I_{m}$.

By Lemma 2.3, we have
\[
\begin{array}{crl}
 -P^{-1}MS^{-1}=-\left(
  \begin{array}{cccccccccccccccc}
   \frac{1}{2+t_{1}}& 0 & 0&...&0\\
  0&\frac{1}{2+t_{2}}&0&...&0\\
  0&0& ...&...&0\\
   0&0& 0&...&\frac{1}{2+t_{m}}\\
 \end{array}
  \right)\left(
  \begin{array}{cccccccccccccccc}
   1^{T}_{t_{1}}& 0 & 0&...&0\\
  0&1^{T}_{t_{2}}&0&...&0\\
  0&0& ...&...&0\\
   0&0& 0&...&1^{T}_{t_{m}}\\
 \end{array}
  \right)\\
  \left(
  \begin{array}{cccccccccccccccc}
   (L_{H_{1}}+I_{t_{1}}-\frac{1}{2+t_{1}}J_{t_{1}})^{-1}& 0 & 0&...&0\\
  0& (L_{H_{2}}+I_{t_{2}}-\frac{1}{2+t_{2}}J_{t_{2}})^{-1}&0&...&0\\
  0&0& ...&...&0\\
   0&0& 0&...& (L_{H_{m}}+I_{t_{m}}-\frac{1}{2+t_{m}}J_{t_{m}})^{-1}\\
 \end{array}
  \right)
\end{array}
\]
\begin{eqnarray}
  =-\left(
  \begin{array}{cccccccccccccccc}
   \frac{1}{2}1^{T}_{t_{1}}& 0 & 0&...&0\\
  0&\frac{1}{2}1^{T}_{t_{2}}&0&...&0\\
  0&0& ...&...&0\\
   0&0& 0&...& \frac{1}{2}1^{T}_{t_{m}}\\
 \end{array}
  \right)
    =F,
\end{eqnarray}

Similarly, $-S^{-1}M^{T}P^{-1}=N^{T}$,
so
$D^{-1}=
 \left(
  \begin{array}{cccccccccccccccc}
   \frac{1}{2}I_{m}&F\\
   F^{T}&S^{-1}\\
 \end{array}
  \right).$

Next we begin with the computation of $\{1\}$-inverse of $L_{R(G)\ominus\wedge_{i=1}^{m}H_{i}}$.

By Lemma 2.8, we have
\[
\begin{array}{crl}
 H&=&L_{G}+D_{G}-\left(
  \begin{array}{cccccccccccccccc}
  R(G)& 0 \\
 \end{array}
  \right)
 \left(
  \begin{array}{cccccccccccccccc}
   \frac{1}{2}I_{m}&F\\
   F^{T}&T^{-1}\\
 \end{array}
  \right)
  \left(
  \begin{array}{cccccccccccccccc}
    R^{T}(G) \\
  0\\
 \end{array}
  \right)\\
  &=&L_{G}+D_{G}-\left(
  \begin{array}{cccccccccccccccc}
    \frac{1}{2}R(G)&R(G)F\\
 \end{array}
  \right)
  \left(
  \begin{array}{cccccccccccccccc}
     R^{T}(G) \\
  0\\
 \end{array}
  \right)\\
  &=&L_{G}+D_{G}-\frac{1}{2}(D_{G}+A_{G})\\
 &=&\frac{3}{2}L_{G},
\end{array}
\]
so $H^{\#}=\frac{2}{3}L^{\#}_{G}$.

According to Lemma 2.8, we calculate $-H^{\#}BD^{-1}$ and $-D^{-1}B^{T}H^{\#}$.
\[
\begin{array}{crl}
-H^{\#}BD^{-1}&=&-\frac{2}{3}L^{\#}_{G}\left(
  \begin{array}{cccccccccccccccc}
 -R(G)& 0 \\
 \end{array}
  \right)
  \left(
  \begin{array}{cccccccccccccccc}
   \frac{1}{2}I_{m}&F\\
   F^{T}&S^{-1}\\
 \end{array}
  \right)\\
  &=&-\frac{2}{3}L^{\#}_{G}\left(
  \begin{array}{cccccccccccccccc}
 -\frac{1}{2}R(G)&-R(G)F\\
 \end{array}
  \right)
   =\left(
  \begin{array}{cccccccccccccccc}
 \frac{1}{3}L^{\#}_{G}R(G)&\frac{2}{3}L^{\#}_{G}R(G)F\\
 \end{array}
  \right)\\
\end{array}
\]
and
\[
\begin{array}{crl}
-D^{-1}B^{T}H^{\#}=-(H^{\#}BD^{-1})^{T}
  =\left(
  \begin{array}{cccccccccccccccc}
   \frac{1}{3}R^{T}(G)L^{\#}_{G}\\
   \frac{2}{3}F^{T}R^{T}(G)L^{\#}_{G} \\
 \end{array}
  \right).\\
\end{array}
\]
We are ready to compute the $D^{-1}B^{T}H^{\#}BD^{-1}$.
\[
\begin{array}{crl}
D^{-1}B^{T}H^{\#}BD^{-1}
&=&\begin{array}{crl}
 \frac{2}{3}\left(
  \begin{array}{cccccccccccccccc}
     \frac{1}{2}I_{m}&F\\
   F^{T}&S^{-1}\\
 \end{array}
 \right) \left(
  \begin{array}{cccccccccccccccc}
 -R^{T}(G)\\
  0\\
 \end{array}
 \right)L^{\#}_{G}
 \left(
  \begin{array}{cccccccccccccccc}
   -R(G)& 0\\
 \end{array}
 \right)
 \left(
  \begin{array}{cccccccccccccccc}
      \frac{1}{2}I_{m}&F\\
   F^{T}&S^{-1}\\
 \end{array}
 \right)
\end{array}\\
&=&\left(
\begin{array}{cccccccccccccccc}
     \frac{1}{6}R^{T}(G)L^{\#}_{G}R(G)&\frac{1}{3}R^{T}(G)L^{\#}_{G}R(G)F\\
   \frac{1}{3}F^{T}R^{T}L^{\#}_{G}R(G)&\frac{2}{3}F^{T}R^{T}L^{\#}_{G}R(G)F\\
 \end{array}
 \right).
\end{array}\\
\]
Based on Lemma 2.3 and 2.8, the following matrix
\begin{eqnarray}
N=\left(
  \begin{array}{cccccccccccccccc}
  \frac{2}{3}L^{\#}_{G}&\frac{1}{3}L^{\#}_{G}R(G)&\frac{2}{3}L^{\#}_{G}R(G)F\\
   \frac{1}{3}R^{T}(G)L^{\#}_{G}&
   \frac{1}{2}I_{m}+ \frac{1}{6}R^{T}(G)L^{\#}_{G}R(G)&F+\frac{1}{3}R^{T}(G)L^{\#}_{G}R(G)F\\
   \frac{2}{3}F^{T}R^{T}(G)L^{\#}_{G} &F^{T}+\frac{1}{3}F^{T}R^{T}(G)L^{\#}_{G}R(G)&
   S^{-1}+\frac{2}{3}F^{T}R^{T}(G)L^{\#}_{G}R(G)F\\
 \end{array}
  \right)
\end{eqnarray}
is a symmetric $\{1\}$- inverse of $L_{R(G)\ominus\wedge_{i=1}^{m}H_{i}}$.

For any $i,j\in V(G)$, by Lemma 2.1 and the Equation $(4.8)$, we have
\begin{eqnarray*}
r_{ij}(L_{R(G)\ominus\wedge_{i=1}^{m}H_{i}})&=&\frac{2}{3}(L^{\#}_{G})_{ii}+\frac{2}{3}(L^{\#}_{G})_{jj}
-\frac{4}{3}(L^{\#}_{G})_{ij}=\frac{2}{3}r_{ij}(G)
\end{eqnarray*}
as stated in $(i)$.

For any $i,j\in V(H_{k})$$(k=1,2,...,m)$, by Lemma 2.1 and the Equation $(4.8)$, we have
\begin{eqnarray*}
r_{ij}(L_{R(G)\ominus\wedge_{i=1}^{m}H_{i}})&=&(L_{H_{k}}+I_{t_{k}}
-\frac{1}{2+t_{k}}j_{t_{k}})^{-1}_{ii}+
(L_{H_{k}}+I_{t_{k}}-\frac{1}{2+t_{k}}j_{t_{k}})^{-1}_{jj}\\&&
-2(L_{H_{1}}+I_{t_{k}}-\frac{1}{2+t_{k}}j_{t_{k}})^{-1}_{ij}.
\end{eqnarray*}

From the left side of above equation, we can obviously have
\begin{eqnarray*}
r_{ij}(F_{k})&=&((L_{H_{k}}+I_{t_{l}})^{-1})_{ii}+
((L_{H_{k}}+I_{t_{l}})^{-1})_{jj}
-2((L_{H_{k}}+I_{t_{l}})^{-1})_{ij},
\end{eqnarray*}
where $F_{k}=H_{k}\vee \{v\}$, i.e, $F_{k}$
is the graph obtained by
adding new edges from an isolated vetrtex
$v$ to every vertex of $H_{k}$.

For any $i,j\in R(G)$, by Lemma 2.1 and the Equation $(4.8)$,
we have
\begin{eqnarray*}
r_{ij}(R(G)\ominus\wedge_{i=1}^{m}H_{i})
&=&r_{ij}(R(G)).
\end{eqnarray*}
By Lemma 3.1 in $\cite{DC}$, $r_{ij}(R(G))
=\frac{2}{3}r_{ij}(G)$, so $r_{ij}(R(G)\ominus\wedge_{i=1}^{m}H_{i})
=\frac{2}{3}r_{ij}(G)$.

For any $i\in V(G)$, $j\in V(H_{k})$$(k=1,2,...,m)$, since
$i$ and $j$ belong to different components, then by Lemma 2.9,
we have
\begin{eqnarray*}
r_{ij}(R(G)\ominus\wedge_{i=1}^{m}H_{i})
&=&r_{ik}(R(G))+r_{kj}(F_{k}).
\end{eqnarray*}

For any $i\in V(H_{k})$, $j\in V(H_{l})$,  by Lemma 2.9,
we have
\begin{eqnarray*}
r_{ij}(R(G)\ominus\wedge_{i=1}^{m}H_{i})
&=&r_{kl}(R(G))+r_{ik}(F_{k})+r_{jl}(F_{l}).
\end{eqnarray*}

By Lemma 2.4, we have
\begin{eqnarray*}
Kf(L_{R(G)\ominus\wedge_{i=1}^{m}H_{i}})
&=&(n+m+\sum_{i=1}^{m}t_{i})tr( N)-{\bf{1}^{T}}N{\bf{1}}\\
&=&(n+m+\sum_{i=1}^{m}t_{i})\left(\frac{2}{3}tr(L^{\#}_{G})
+tr\left(\frac{1}{2}I_{m}+\frac{1}{6}R^{T}(G)L^{\#}_{G}R(G)\right)+\right.\\
&&\left.+
tr(S^{-1}+\frac{2}{3}F^{T}R^{T}L^{\#}_{G}R(G)F)\right)-
{\bf{1}^{T}}N{\bf{1}}\\
&=&(n+m+\sum_{i=1}^{m}t_{i})\left(\frac{2}{3n}Kf(G)
+\frac{m}{2}+\frac{1}{6}\sum_{i<j,i,j\in E(G)}[(L^{\#}_{G})_{ii}+(L^{\#}_{G})_{jj}\right.\\
&&\left.+2(L^{\#}_{G})_{ij}]+tr\left(S^{-1}+\frac{2}{3}F^{T}R^{T}(G)L^{\#}_{G}R(G)F\right)\right)
-{\bf{1}^{T}}N{\bf{1}}
\end{eqnarray*}
By Lemma 2.5, we get
\begin{eqnarray*}
Kf(L_{R(G)\ominus\wedge_{i=1}^{m}H_{i}})
&=&(n+m+\sum_{i=1}^{m}t_{i})\left(\frac{2}{3n}Kf(G)
+\frac{m}{2}+\frac{1}{6}\sum_{i<j,i,j\in E(G)}[2(L^{\#}_{G})_{ii}+2(L^{\#}_{G})_{jj}\right.\\
&&\left.-r_{ij}(G)]+tr\left(S^{-1}+\frac{2}{3}F^{T}R(G)^{T}L^{\#}_{G}R(G)F\right)\right)
-{\bf{1}^{T}}N{\bf{1}}\\
&=&(n+m+\sum_{i=1}^{m}t_{i})\left(\frac{2}{3n}Kf(G)
+\frac{m}{2}+\frac{1}{3}tr(D_{G}L^{\#}_{G})-\frac{n-1}{6}\right.\\
&&\left.+tr\left(S^{-1}+\frac{2}{3}F^{T}R(G)^{T}L^{\#}_{G}R(G)F\right)\right)
-{\bf{1}^{T}}N{\bf{1}}
\end{eqnarray*}
Note that the eigenvalues of $(L_{H_{i}}+I_{t_{i}}-\frac{1}{2+t_{i}}j_{t_{i}})$
$(i=1,2,...,m)$ are $\mu_{1}(H_{i})+1, \mu_{2}(H_{i})+1,...,\mu_{t_{i}}(H_{i})+1$. Then
\begin{eqnarray}
tr(S^{-1})=
\sum_{i=1}^{m}\sum_{j=1}^{t_{i}}\frac{1}{\mu_{i}(H_{j})+1}.
\end{eqnarray}
By Lemma 2.2, $L_{G}^{\#}\bf{1}$$=0$ and $({\bf{1}^{T}}\left(R^{T}(G)L^{\#}_{G}Q\right){\bf{1}})^{T}
={\bf{1}^{T}}\left(Q^{T}L^{\#}_{G}R(G)\right){\bf{1}}$, then
\begin{eqnarray*}
{\bf{1}^{T}}N{\bf{1}}&=&
\frac{m}{2}+
\frac{1}{6}{\bf{1}^{T}}\left(R^{T}(G)L^{\#}_{G}R(G)\right){\bf{1}}+
{\bf{1}^{T}} F{\bf{1}}+{\bf{1}^{T}}F^{T}{\bf{1}}
\\
&&+\frac{2}{3}{\bf{1}^{T}} R^{T}(G)L^{\#}_{G}R(G)F{\bf{1}}+{\bf{1}^{T}}S^{-1}{\bf{1}}
+\frac{2}{3}{\bf{1}^{T}}\left(F^{T}R^{T}(G)L^{\#}_{G}R(G)F\right){\bf{1}}.
\end{eqnarray*}
Note that $R(G)\bf{1}=\pi$, where $\pi^{T}=(d_{1},d_{2},...,d_{n})$,
then ${\bf{1}^{T}}\left(R^{T}(G)L^{\#}_{G}R(G)\right){\bf{1}}
=\pi^{T}L^{\#}_{G}\pi$, so
\begin{eqnarray}
{\bf{1}^{T}}N{\bf{1}}&=&\frac{m}{2}+\frac{1}{6}\pi^{T}L^{\#}_{G}\pi+
\pi^{T}L^{\#}_{G}Q(G){\bf{1}}+{\bf{1}^{T}}T^{-1}{\bf{1}}+{\bf{1}^{T}}
\left(Q^{T}L^{\#}_{G}Q\right){\bf{1}}.
\end{eqnarray}
Let $R_{i}=L_{H_{i}}+I_{t_{i}}-\frac{1}{2+t_{i}}j_{t_{i}}$$(i=1,2,...,m)$, then
\begin{eqnarray*}
{\bf{1}^{T}}S^{-1}{\bf{1}^{T}}&=&\left(
  \begin{array}{cccccccccccccccc}
 {\bf{1}^{T}_{t_{1}}}&{\bf{1}^{T}_{t_{2}}}&\cdots &{\bf{1}^{T}_{t_{m}}}\\
 \end{array}
  \right)
  \left(
  \begin{array}{cccccccccccccccc}
R_{1}^{-1}& 0 & 0&...&0\\
  0&R_{2}^{-1}&0&...&0\\
  0&0& ...&...&0\\
   0&0& 0&...&R_{m}^{-1}\\
 \end{array}
  \right)\left(
  \begin{array}{cccccccccccccccc}
   {\bf{1}_{t_{1}}}\\
  {\bf{1}_{t_{2}}}\\
  \cdots\\
  {\bf{1}_{t_{m}}}\\
 \end{array}
  \right)\\
  \end{eqnarray*}
  \begin{eqnarray}
 &=&\sum_{i=1}^{m}{\bf{1}^{T}_{t_{i}}}(L_{H_{i}}+I_{t_{i}}-\frac{1}{2+t_{i}}j_{t_{i}})^{-1}{\bf{1}_{t_{i}}}
 =\frac{1}{2}\sum_{i=1}^{m}t_{i}(2+t_{i}),
\end{eqnarray}
and
\begin{eqnarray*}
{\bf{1}^{T}}F^{T}&=&\frac{1}{2}\left(
  \begin{array}{cccccccccccccccc}
 {\bf{1}^{T}_{t_{1}}}&{\bf{1}^{T}_{t_{2}}}&\cdots &{\bf{1}^{T}_{t_{m}}}\\
 \end{array}
  \right)
  \left(
  \begin{array}{cccccccccccccccc}
 {\bf{1}_{t_{1}}}& 0 & 0&...&0\\
  0& {\bf{1}_{t_{2}}}&0&...&0\\
  0&0& ...&...&0\\
   0&0& 0&...& {\bf{1}_{t_{m}}}\\
 \end{array}
  \right)\\
    \end{eqnarray*}
\begin{eqnarray}
 &=&\frac{1}{2}(t_{1},0,...,0,t_{2},0,...,0...,t_{m})=\frac{1}{2}\delta^{T}.
\end{eqnarray}

Plugging $(4.9),(4.10),(4.11)$ and $(4.12)$ into $Kf(L_{R(G)\ominus\wedge_{i=1}^{m}H_{i}})$,
we obtain the required result in $vi)$.

\section{Conclusion}
In this paper,
using the Laplacian generalized inverse
approach, we obtained the resistance distance
and Kirchhoff indices of $R(G)\boxdot \wedge_{i=1}^{n}H_{i}$
and $R(G)\ominus\wedge_{i=1}^{m}H_{i}$ whenever $G$ and
$H_{i}$ are arbitrary graph. These results
generalize the existing results in $\cite{HPZ}$.

This article has been reviewed in Filomat on November 19, 2017.

\vskip 0.1in
\noindent{\bf Acknowledgment:}
This work was supported by the National Natural Science Foundation of China (No.11461020) and the Research Foundation of the Higher Education Institutions of Gansu Province, China (2018A-093).

\end{document}